# Another Look at Random Infinite Divisibility

## S. Satheesh

Department of Statistics
Cochin University of Science and Technology
Cochin 682 022, India.
ssatheesh@sancharnet.in

**Abstract:** Surveying random ($N$) infinite divisibility, where $N$ is a non-negative integer-valued random variable, one comes across a definition where the class of probability generating functions $\{P_\theta\}$ of $N$ forms a commutative semi-group. We show that this assumption is not natural since it captures the notion of stable rather than infinitely divisible laws. Also it rules out any $N$ having an atom at the origin. In this paper we propose another description of N-infinite divisibility, which is applicable to a wider class of $N$. This is $\varphi$-infinite divisibility where $\varphi$ is a Laplace transform. Here the commutative semi-group assumption is not required and includes $N$ having an atom at zero. We discuss attraction and partial attraction in this et up. The role of the divisibility of $N$ and $\varphi$ on that of the $\varphi$-ID law is explored and we give a method to construct class-L laws.

**Keywords and Phrases:** infinite divisibility, stability, semi-stability, attraction, partial attraction, geometric infinite divisibility, random sum, N-infinite divisibility, N-attraction, de-Finetti theorem, characteristic function, Laplace transform, probability generating function, transfer theorem.

## 1. Introduction

In the classical summation scheme a characteristic function (CF) $f(t)$ is infinitely divisible (ID) if for every $n \geq 1$ integer there exists a CF $f_n(t)$ such that

$$f(t) = \{f_n(t)\}^n. \tag{1.1}$$

The classical de-Finetti theorem for ID laws states that $f(t)$ is ID if and only if

$$f(t) = \underset{n \to \infty}{Lt}\ exp\{-a_n(1-h_n(t))\} \tag{1.2}$$

where $\{a_n\}$ are some positive constants and $h_n(t)$ are CFs.

## Another Look at Random Infinite Divisibility

Klebanov, *et al.* (1984) have extended the notion of ID laws to geometric (with mean $1/p$) summation schemes introducing geometrically infinitely divisible (GID) laws. According to this $f(t)$ is GID if for every $p \in (0,1)$ there exists a CF $f_p(t)$ such that

$$f(t) = \sum_{n=1}^{\infty} (f_p(t))^n \, p(1-p)^{n-1} \qquad (1.3)$$

the geometric law being independent of the distribution of $f_p(t)$. They also proved an analogue of the de-Finetti theorem in the context, *viz.* $f(t)$ is GID if and only if

$$f(t) = \operatorname*{Lt}_{n \to \infty} 1/\{1 + a_n(1 - h_n(t))\}, \qquad (1.4)$$

where $\{a_n\}$ and $\{h_n(t)\}$ are as in (1.2). Consequently $f(t)$ is GID if and only if

$$f(t) = 1/\{1 - \log \omega(t)\}, \qquad (1.5)$$

where $\omega(t)$ is a CF that is ID. Subsequently Sandhya (1991) (also reported in Sandhya and Pillai (1999)), Mohan, *et al.* (1993), Ramachandran (1997), Gnedenko and Korolev (1996) and Klebanov and Rachev (1996) have discussed attraction and the first three works that of partial attraction in geometric sums.

Sandhya (1991, 1996), Gnedenko and Korolev (1996), Klebanov and Rachev (1996) and Bunge (1996) have extended the notion of infinite divisibility to random sum (N-sum) schemes. Sandhya (1991, 1996) defined N-ID laws as: $f(t)$ is N-ID, where $N$ is a positive integer valued r.v with probability generating function (PGF) $P_\theta$, if there exists a CF $f_\theta(t)$ such that

$$f(t) = P_\theta\{f_\theta(t)\} \quad \text{for every } \theta \in \Theta. \qquad (1.6)$$

Here $N$ is assumed to have a finite mean and is independent of the distribution of $f_\theta(t)$ and $\Theta$ is the parameter space of $\theta$. She also noticed that when $f(t)$ and $f_\theta(t)$ are of the same type (*ie.* when $f(t)$ is N-sum-stable) (1.6) is an Abel equation and discussed N-semi-stable laws as well. One drawback of this definition is that it is not constructive. However she did give two examples of non-geometric laws for $N$.

Gnedenko and Korolev (1996) and Klebanov and Rachev (1996) went further from (1.6), proving the de-Finetti analogue for N-sums *viz.* $f(t)$ is N-ID if and only if

$$f(t) = \underset{n\to\infty}{Lt}\, \varphi\{a_n(1- h_n(t))\} \tag{1.7}$$

where $\varphi$ is a Laplace transform (LT) and $\{a_n\}$ and $\{h_n(t)\}$ are as in (1.2). They then concluded that $f(t)$ is N-ID if and only if

$$f(t) = \varphi\{-\log \omega(t)\} \tag{1.8}$$

where $\omega(t)$ is the CF of an ID law. They first developed N-normal laws and then proceeded to define N-ID laws using the Poincare equation (Abel equation). This approach required that the family of PGFs $\{P_\theta:\theta\in\Theta\}$ formed a commutative semi-group with respect to the operation of convolution in addition to the existence of the mean of $N$. Here $P_\theta$ and $\varphi$ (in (1.6), (1.7) and (1.8)) are related by

$$P_\theta(z) = \varphi\{\varphi^{-1}(z)/\theta\},\ \theta\in\Theta. \tag{1.9}$$

Bunge (1996) also arrived at (1.8) under the same assumptions but the arguments were based on Levy processes instead of proving the de-Finetti analogue enroute. The relation (1.9) holds true in this development also. But Gnedenko and Korolev (1996) in addition, have explicitly proved that N-ID laws and they alone can be the limit laws of N-sums by invoking a transfer theorem. Gnedenko and Korolev (1996) and Klebanov and Rachev (1996) also discuss attraction in this N-sum scheme.

All these authors also discuss N-stable laws while Klebanov and Rachev (1996) N-semi-stable laws as well, the requirement being that $\omega(t)$ in (1.8) must be the CF of a stable or a semi-stable law. For N-stable laws Kozubowski and Panorska (1996) have reached the same conclusion under the assumption that $N_\theta \xrightarrow{p} \infty$ as $\theta\downarrow 0$, without requiring that the PGF of $N_\theta$, $\{P_\theta:\theta\in\Theta\}$ form a commutative semi-group.

# Another Look at Random Infinite Divisibility

This paper is motivated by the following observations. The development of N-ID laws based on the requirement that $\{P_\theta : \theta \in \Theta\}$ formed a commutative semi-group does not appear to be natural (see Remark.2.1) and is rather restrictive. From a methodological perspective also this is important as N-ID laws approximate N-sums and one cannot expect the distribution of *N* (sample size) to satisfy such an assumption. In fact the assumption rules out any $P_\theta$ having an atom at the origin. A case in point is when the distribution of *N* is negative binomial and we do not have a corresponding N-ID law. On the other hand, when $\varphi$ is gamma one expects the corresponding φ-ID law to approximate a negative binomial sum. But in this case $P_\theta$ corresponds to a Harris law see example.2.4. Another consequence of this assumption is that though the classical de-Finetti theorem describes an N-sum it does not fit in to this scheme of development while its analogues are the key results in the development of GID and N-ID laws.

Generalizing different descriptions of ID laws has done generalization of ID laws to N-sums. Those surveyed here are generalizations of (1.1) while generalizations of triangular null array scheme are also available in the literature; see references in Klebanov and Rachev (1996). Here we generalize the fact that the class of ID laws coincides with the class of limit laws of compound Poisson laws (1.2). This generalization straight away describes an N-sum and we name them as φ-ID laws. Here the assumptions (i) finite mean for *N* (ii) $\{P_\theta : \theta \in \Theta\}$ form a commutative semi-group and (iii) $N_\theta \xrightarrow{p} \infty$ as $\theta \downarrow 0$, are done away with and thus the requirements are much weaker. Arriving at (1.8) we show that the known results are true in this setup also. We then show that our description of φ-ID laws lends itself to a much wider class of PGFs for *N*. However, we will show (Remark.2.1) that when we consider N-sum-stable laws as in (1.6) it is natural that $\{P_\theta : \theta \in \Theta\}$ formed a commutative semi-group. These are done in Section.2. Certain divisibility properties of φ-ID laws are proved in Section.3. We develop and discuss φ-attraction and partial φ-attraction in Section.4. This development of φ-ID laws is motivated by Feller's proof of Bernstein's theorem for completely monotone functions (Feller, 1971, p.440) and the (transfer) Theorem.4.1.2 in Gnedenko and Korolev (1996, p.98) (hereafter referred to as GK (1996)).

Subsequently, Satheesh (2002) studied randomization and mixtures of ID and max-ID laws and related processes, Satheesh and Sandhya (2003) more on φ-ID and φ-max-ID laws, Satheesh (2003) operator φ-stable laws and Satheesh and Sandhya (2004) operator φ-semi-stable and max-self-decomposable laws.

## 2. φ-ID Laws

We begin with the observation that perhaps the commutative semi-group assumption is not a natural setting for the description of N-ID laws.

**Remark.2.1** Suppose that the $N_\theta$-sum of $\varphi(s)$ is of the same type. That is $\varphi(s) = P_\theta\{\varphi(\theta s)\}$ for every $\theta \in \Theta$. Equivalently, $\varphi(s/\theta) = P_\theta\{\varphi(s)\}$, implying $P_\theta(z) = \varphi\{\varphi^{-1}(z)/\theta\}$, $\theta \in \Theta$ when $s = \varphi^{-1}(z)$. Thus (1.9) captures the structure of N-sum-stability. By corollary.4.6.1 in GK (1996, p.141) (1.9) is equivalent to the assumption that $\{P_\theta(z):\theta \in \Theta\}$ form a commutative semi-group. Clearly the converse also holds. The argument is also true for an N-stable CF of the form (1.8) where $\omega(t)$ is the CF of a stable law. Thus the development of N-ID laws from N-normal laws, based on the commutative semi-group assumption captures the fact that normal laws are stable and not just ID. This is why this assumption got in to the scheme of N-ID laws as an extension of N-normal laws.

**Note.** An implication of Theorem.1a in Feller (1971, p.432) is: a LT has no real zeroes; thus:

**Remark.2.2** Another implication of (1.9) is that $N$ cannot have an atom at the origin. This is because at $z = 0$, the RHS of (1.9) is zero and hence $P_\theta(o) = 0$. Hence the conclusion. In particular we cannot have N-sum stability when $N$ is Poisson or negative binomial.

The following lemma was used to define discrete analogue of the notion of distributions of the same type in Satheesh and Nair (2002) and Satheesh, *et al.* (2002).

Another Look at Random Infinite Divisibility

**Lemma.2.1** If $\varphi(s)$, $s>0$ is a LT then $P_\theta(s) = \varphi\{(1-s)/\theta\}$, $0<s<1$, $\theta>0$ is a PGF.

**Definition.2.1** Let $\varphi$ be a LT. A CF $f(t)$ is $\varphi$-ID if

$$f(t) = \underset{n\to\infty}{Lt}\ \varphi\{a_n(1-f_n(t))\} \qquad (2.1)$$

where $\{a_n\}$ are some positive constants and $\{f_n(t)\}$ are CFs and the distributions of $\{f_n(t)\}$ and $\varphi\{a_n(1-s)\}$ are independent for each $n$.

Since for each $n$, $\varphi\{a_n(1-f_n(t))\}$ represent an N-sum of r.vs with CF $f_n(t)$, where the PGF of $N$ is $\varphi\{a_n(1-s)\}$, the definition identifies $f(t)$ as the limit of a sequence of N-sums. Notice that this requirement is weaker than that $f(t)$ must be expressible as an $N_\theta$-sum of i.i.d r.vs for each $\theta \in \Theta$. Here the distributions of the r.v $N$ and the component r.vs (with CF $f_n(t)$) are explicit unlike those in GK (1996), Klebanov and Rachev (1996) or Bunge (1996). The following is analogous to (1.8).

**Theorem.2.1:** A CF $f(t)$ is $\varphi$-ID if and only if $f(t) = \varphi\{\psi(t)\}$ and $\omega(t) = e^{-\psi(t)}$ is a CF that is ID.

**Proof:** We have $f(t) = \underset{n\to\infty}{Lt}\ \varphi\{a_n(1-f_n(t))\} = \varphi(-\log\ \underset{n\to\infty}{Lt}\ exp(-\{a_n(1-f_n(t))\}))$,

as in GK (1996, p.147). Now invoking the classical de-Finetti theorem we have a CF $\omega(t) = e^{-\psi(t)}$ that is ID and,

$$f(t) = \varphi\{-\log\omega(t)\} = \varphi\{\psi(t)\} . \qquad (2.2)$$

The converse now follows by invoking the same. The proof is thus complete.

**Theorem.2.2:** A CF $f(t)$ which is the limit of a sequence of CFs that are $\varphi$-ID is again $\varphi$-ID.

S. SatheeshS. Satheesh

**Proof.** Let $\{g_n(t)\}$ be a sequence of CFs that are $\varphi$-ID and converging to the CF $f(t)$. Then there exists CFs $\{\omega_n(t)\}$ that are ID, and representations corresponding to (2.2) for each $\{g_n(t)\}$. Hence,

$$f(t) = \underset{n\to\infty}{Lt}\ g_n(t) = \underset{n\to\infty}{Lt}\ \varphi\{-\log \omega_n(t)\} = \varphi\{-\log \underset{n\to\infty}{Lt}\ \omega_n(t)\} = \varphi\{-\log \omega(t)\}.$$

Now, $\omega(t)$ is ID being the limit of a sequence of CFs that are ID. Hence $f(t)$ is $\varphi$-ID.

**Theorem.2.3:** For a CF $g(t)$ and $a>0$, $\varphi\{a(1-g(t))\}$ is $\varphi$-ID.

**Proof:** Follows since $exp\{-a(1-g(t))\}$ is ID.

**Theorem.2.4:** A CF $f(t)$ that is $\varphi$-ID has no real zeroes.

**Proof:** Suppose that $f(t)$ has a real zero, say $\tau$. Then, $0 = f(\tau) = \varphi\{-\log \omega(\tau)\}$ where $\omega(.)$ is a CF that is ID. Hence, $-\log \omega(\tau) = \varphi^{-1}(o) = \infty$ (see the note after Remark.2.1). This implies

$$\omega(\tau) = exp(-\infty) = 0.$$

That is $\omega(.)$, a CF that is ID, also has a real zero that is a contradiction. Hence the proof.

**Lemma.2.2:** Given a r.v $U$ with LT $\varphi$ there exists an integer valued r.v $N_\theta$, such that

$$\theta N_\theta \xrightarrow{L} U \text{ as } \theta \to 0.$$

**Proof:** By Lemma.2.1, $P_\theta(s) = \varphi\{(1-s)/\theta\}$, $0<s<1$, $\theta>0$ is a PGF. Let $N_\theta$ be the corresponding r.v. The LT of $\theta N_\theta$ is

$$\varphi\{(1-e^{-v\theta})/\theta\},\ v>0.$$

129129

Now as given in Feller (1971, p.440)

$$\underset{\theta \to 0}{Lt}\ \varphi\{(1- e^{-v\theta})/\theta\} = \varphi(v),$$

completing the assertion.

**Theorem.2.5:** $Z$ is the limit law, as $\theta \to 0$, of $N_\theta$-sums of i.i.d r.vs (whose classical-sum has a limit law) where the PGF of $N_\theta$ is $\varphi\{(1-s)/\theta\}$, if and only if $Z$ is φ-ID.

**Proof:** The assertion follows (along the lines as in GK (1996, p.148)) from Lemmata.2.1, 2.2, Theorem.2.1, and the transfer theorem in GK (1996, p.98).

**Example.2.1:** From definition.2.1 one sees that the classical de-Finetti theorem readily follows when $\varphi$ is a degenerate LT (as the limit of Poisson sums).

**Example.2.2:** GID laws are φ-ID when $\varphi$ is exponential, as the limit of geometric (with mean $1/p$) sums. Now the geometric law can also have the support $\{0,1,2,\ldots\}$ and one may identify $\theta = p/q$, $q = 1-p$, $p$ being the parameter of the geometric law in this case.

**Example.2.3:** When $\varphi$ is gamma the φ-ID law approximate negative binomial sums.

Now from Theorem.2.1 we define φ-stable and φ-semi-stable laws (See also Klebanov and Rachev (1996)) and then describe them in the spirit of Definition.2.1.

**Definition.2.2:** A CF $f(t)$ is φ-stable (φ-semi-stable) if and only if $f(t) = \varphi\{-\log \omega(t)\}$ where $\omega(t)$ is a stable (semi-stable) CF. The following lemma is now clear.

**Lemma.2.3:** A CF $f(t)$ is φ-ID (respectively φ-stable or φ-semi-stable) if and only if $exp\{-\varphi^{-1}(f(t))\}$ is ID (respectively stable or semi-stable).

**Theorem.2.6:** A CF $f(t)$ is $\varphi$-stable if and only if

$$f(t) = \underset{n \to \infty}{Lt}\ \varphi\{n\,(1 - h(t/a_n)) + it\mu_n\}$$

for some sequences of real constants $\{a_n > 0\}$, $\{\mu_n\}$ and a CF $h(t)$.

**Proof:** The assertion follows essentially from the fact that

$$\underset{n \to \infty}{Lt}\ \{n\,(1 - h(t/a_n)) + it\mu_n\} = \psi(t)$$

implies $e^{-\psi(t)}$ is stable and the converse is also true; see Ramachandran (1997). This result generalizes his Proposition.2.1.

**Corollary.2.1:** In this theorem $\mu_n = 0\ \forall\,n$, if and only if $f(t)$ is $\varphi$-strictly stable. The conclusion follows from Theorem.2.2 in Mohan, *et al.* (1993).

**Theorem.2.7:** A CF $f(t)$ is $\varphi$-semi-stable if and only if

$$f(t) = \underset{k \to \infty}{Lt}\ \varphi\{n_k\,(1 - h(t/a_k))\}$$

where $\{a_k\}$ are some positive constants, $\{n_k\}$ a sequence of positive integers satisfying $\underset{k \to \infty}{Lt}\ \inf\{n_k/n_{k+1}\} > 0$ and $h(t)$ a CF.

**Proof:** The assertion follows from the fact that under the specified conditions $\underset{k \to \infty}{Lt}\ \{n_k\,(1 - h(t/a_k))\} = \psi(t)$ imply that $e^{-\psi(t)}$ is semi-stable and the converse is also true. See Theorem.3.2 in Mohan, *et al.* (1993).

**Corollary.2.2:** If in this theorem $h(t) = f(t)$ and $\underset{k \to \infty}{Lt}\ \{n_k/n_{k+1}\} = 1$ then $f(t)$ is $\varphi$-strictly stable. The conclusion follows from Theorem.2.3 in Mohan, *et al.* (1993).

Theorem.2.5 enables one to envisage φ-ID laws as approximating N-sums when the sample size distribution *N* has PGF $\varphi\{(1-s)/\theta\}$ as $\theta \to 0$. However, here $P_\theta(o)>0$ and so distributions having no atom at the origin in general are excluded and hence in particular the GID laws. Notice that the class of PGFs $P_\theta(z) = \varphi\{\varphi^{-1}(z)/\theta\}$, $\theta \in \Theta$, under the commutative semi-group assumption have no atom at zero (Remark.2.2). In order to widen the applicability of our approach we need bring in more PGFs to the folder of Theorem.2.5. In Theorem.2.8 we achieve this.

**Lemma.2.4:** $\wp_\varphi = \{P_\theta(s) = s^j \varphi\{(1-s^m)/\theta\}, 0<s<1, j \geq 0 \text{ \& } m \geq 1 \text{ integers and } \theta > 0\}$ describes a class of PGFs for any given LT $\varphi$.

**Proof:** Follows essentially from Lemma.2.1 by noticing that the requirements of absolute monotonicity and of $P_\theta(1) = 1$ are intact.

**Lemma.2.5:** Given a r.v $U$ with LT $\varphi$, the integer valued r.vs $N_\theta$ with PGF $P_\theta$ in the class $\wp_\varphi$ described in Lemma.2.4 satisfy

$$\theta N_\theta \xrightarrow{L} mU \text{ as } \theta \to 0.$$

**Proof:** Now the LT of $\theta N_\theta$ is

$$e^{-vj\theta} \varphi\{(1-e^{-vm\theta})/\theta\}, v>0$$

and

$$\underset{\theta \to 0}{Lt} \{e^{-vj\theta} \varphi((1-e^{-vm\theta})/\theta)\} = \varphi(mv),$$

which now follows along the same lines as that of Lemma.2.2.

Satheesh and Sandhya (2004) have shown that these $N_\theta \xrightarrow{p} \infty$ as $\theta \downarrow 0$.

**Theorem.2.8:** The limit law, as $\theta \to 0$, of $N_\theta$-sums of i.i.d r.vs (whose classical-sum has a limit law), where the PGF of $N_\theta$ is a member of $\wp_\varphi$, is necessarily φ-ID. Conversely, each φ-ID law can be obtained as the limit law of $N_\theta$-sums of some i.i.d r.vs for every $N_\theta$ whose PGF is in $\wp_\varphi$.

**Proof:** Follows from the transfer theorem, Theorem.2.1, Lemmata.2.4 and 2.5. Notice that if a CF $\omega$ is ID then $\omega^m$, $m>0$ is also ID. Thus:

**Definition.2.3:** A CF $f(t)$ is φ-ID if there exists a CF $h_\theta(t)$, a PGF $P_\theta \in \wp_\varphi$ that is independent of $h_\theta$ for each $\theta \in \Theta$, such that $P_\theta\{h_\theta(t)\} \to f(t)$ as $\theta \to 0$ through a sequence $\{\theta_n\}$.

**Remark.2.3:** The uniqueness of $P_\theta$ given $\varphi$ in (1.9), is not here as in the commutative semi-group approach. However, here the φ-ID law approximates the N-sum when the PGF of $N$ is a member of $\wp_\varphi$. Notice that the de-Finetti analogue for GID laws in (1.4) can also be seen as describing a geometric (on $\{0,1,2, \ldots\}$) sum while the summation description corresponds to a geometric law on $\{1,2, \ldots\}$. In the case of GID laws they coincide under the commutative semi-group assumption (implicit in the proof of (1.4)).

**Example.2.4:** When $N_\theta$ is a Harris $(a,m)$ law with PGF, $P_a(s) = s/\{a - (a-1)s^m\}^{1/m}$, $m>0$ integer, $a>1$, setting $\theta = 1/a$, $\underset{\theta \to o}{Lt}\ \theta N_\theta = U$ and $U$ has a gamma$(m,1/m)$ law. As the solution of the Poincare equation this has also been given in Klebanov and Rachev (1996, p.489). Notice also that when $\varphi$ is gamma$(m,\beta)$, under N-sum stability as in (1.6) $N_\theta$ is Harris$(a,m)$ by (1.9) and $\beta = 1/m$, see Satheesh, *et al.* (2002). (See also Remark.3.2).

### 3. Divisibility Properties of φ-ID Laws

If $\varphi$ is ID then the φ-ID law with CF $f(t) = \varphi\{-\log \omega(t)\}$ in (2.2) is also ID by Property.4.6.3 of GK (1996, p.148). But is it necessary? From Theorem.2.3 we know that $\varphi\{a(1- g(t))\}$ is φ-ID. Alternately, the problem may be posed whether a mixture of compound Poisson laws can be ID while $\varphi$ is not ID. Or still, whether a mixture of Poisson laws can be ID while $\varphi$ is not ID. By Kallenberg's counter example (Grandell, 1997, p.28) a mixture of Poisson laws can be ID while $\varphi$ is not ID. Thus:

**Property.3.1:** A φ-ID law can be ID even when $\varphi$ is not ID.

A similar question can be asked in the context of N-ID laws as described by (1.6). When the PGF $P_\theta$ is ID the N-sum is also ID by Feller (1971, p.464). We know that GID laws are ID. But here *N* is geometric on {1,2, ….} and its PGF is not ID (compound Poisson), Feller (1968, p.290). (See Satheesh (2003) for more on this). Thus:

**Property.3.2:** An N-ID law can be ID even when the PGF of *N* is not ID.

In the case of class-L (*L*) laws we have from (Sandhya and Satheesh, 1996): A CF *f(t)* that is N-sum strictly stable is in *L* if *N* is positive. Though geometrically stable laws are in *L*, the PGF of this geometric law is not in discrete-*L* (Steutel and van Harn, 1979) as it is not even ID. Thus:

**Property.3.3** An N-ID law can be in *L* even when the PGF of *N* is not in discrete-*L*.

Theorem.2.1 of Ramachandran (1997) gives sufficient conditions for a geometric stable law to be in *L*. The next result gives a sufficient condition for φ-ID laws.

**Property.3.4** A CF *f(t)* that is φ-ID is in *L* if it is φ-strictly stable and $\varphi$ is in *L*. Consequently *f(t)* is unimodal and absolutely continuous.

**Proof.** If $\varphi$ is in *L* then there is another LT $\varphi_c$ such that

$$\varphi(s) = \varphi(cs)\, \varphi_c(s),\ s>0\ \text{for every}\ 0<c<1.$$

Now, since *f(t)* is φ-strictly stable

$$f(t) = \varphi(\psi(t)) = \varphi(c\psi(t))\, \varphi_c(\psi(t)),$$

where $\psi(t) = \lambda|t|^\alpha e^{-i\theta \operatorname{sgn}(t)}$, $\lambda>0$, $|\theta| \leq \min(\pi\alpha/2, \pi-\pi\alpha/2)$, $0<\alpha\leq 2$.

This proves the assertion and the consequences follow.

**Remark.3.1:** Since gamma $(\lambda,\nu)$ law is in *L*, the corresponding $\varphi$-strictly stable law is also in *L*. Thus the generalized Linnik (GL $(\alpha,\theta,\nu)$) law of Erdogan and Ostrowski (1997) is in *L* and hence unimodal and absolutely continuous. This argument is much simpler than the one given by them. Membership in *L* or not is also important from a modeling perspective as such distributions can model AR processes, Gaver and Lewis (1980). Property.3.4 also enables one to construct distributions in *L*. This property generalizes Theorems 2.3, 2.4 and 2.5 in Satheesh, *et al.* (2002).

**Remark.3.2:** Since the limit law of $N(a)/a$ as $a \to \infty$ is gamma$(m,1/m)$ when $N(a)$ is Harris$(a,m)$ the corresponding $\varphi$-stable law is the gamma$(m,1/m)$ mixture of stable laws. The stability of the generalized Linnik law in the setup as a negative binomial $(p)$ sum has been mentioned in Kozubowski and Panorska (1998) since in this case also the limit law of $pN(p)$ as $p \to 0$ is gamma$(1,1/m)$. However, a GL$(\alpha,\theta,1/m)$ law is not N-sum stable (as conceived in (1.6)) w.r.t a negative binomial$(p)$ r.v. (See also Example.2.4).

## 4. $\varphi$-attraction and Partial $\varphi$-attraction

In the classical summation scheme a CF $g(t)$ belongs to the domain of attraction (DA) of the CF $\omega(t)$ if there exist sequences of real constants $a_n > 0$ and $b_n$ such that

$$\text{as } n \to \infty, \; \{g(t/a_n) \, exp(-itb_n)\}^n \to \omega(t) \text{ for all } t \in \mathbf{R}.$$

By setting $g_n(t) = g(t/a_n) \, exp(-itb_n)$ this is equivalent to;

$$\text{as } n \to \infty, \; exp\{-n(1 - g_n(t))\} \to exp\{-\psi(t)\} = \omega(t) \; \forall \; t \in \mathbf{R}. \qquad (4.1)$$

When convergence is possible only as *n* runs through a subsequence $\{n_k\}$ of positive integers we have: A CF $g(t)$ belongs to the domain of partial attraction (DPA) of the CF $\omega(t)$ if there exists sequences of real constants $a_k > 0$ and $b_k$ and $n_k \to \infty$ such that

$$\text{as } k \to \infty, \ exp\{-n_k(1 - g_k(t))\} \to exp\{-\psi(t)\} = \omega(t) \ \forall \ t \in \mathbf{R}. \quad (4.2)$$

If $b_n = 0$ we have the corresponding strict-sense domains. We now extend the notions w.r.t the PGFs in Lemma.2.1 as:

**Definition.4.1:** A CF $g(t)$ belongs to the domain of φ-attraction (Dφ-A) of the CF $f(t)$ if there exist a sequence of PGFs $\varphi\{n(1-s)\}$ independent of $g(t)$, and sequences of constants $a_n > 0$ and $b_n \in \mathbf{R}$ such that

$$\text{as } n \to \infty, \ \varphi\{n(1 - g_n(t))\} \to f(t) \text{ for all } t \in \mathbf{R}. \quad (4.3)$$

**Definition.4.2:** A CF $g(t)$ belongs to the domain of partial φ-attraction (DPφ-A) of the CF $f(t)$ if there exist a sequence of PGFs $\varphi\{n_k(1-s)\}$ independent of $g(t)$, and sequences of constants $a_k > 0$ and $b_k \in \mathbf{R}$ and $n_k \to \infty$ such that

$$\text{as } k \to \infty, \ \varphi\{n_k(1 - g_k(t)) \to f(t) \text{ for all } t \in \mathbf{R}. \quad (4.4)$$

We now have the following results. Recall that a φ-ID CF has no real zeroes.

**Theorem.4.1:** A CF $f(t)$ has a nonempty Dφ-A if and only if it is φ-stable.

**Proof:** Suppose $f(t)$ is φ-stable. Then there exists a stable CF $e^{-\psi(t)}$ such that $f(t) = \varphi\{\psi(t)\}$. Since $e^{-\psi(t)}$ is stable it has a DA and let $g(t)$ be a member of it. Hence;

$$\text{as } n \to \infty, \ n(1 - g_n(t)) \to \psi(t).$$

Thus, $\varphi\{n(1 - g_n(t))\} \to \varphi\{\psi(t)\} = f(t)$ implying $g(t)$ is a member of the Dφ-A of $f(t)$.

Conversely, if $g(t)$ belongs to the D$\varphi$-A of $f(t)$ then

$$f(t) = \underset{n\to\infty}{Lt}\ \varphi\{n(1 - g_n(t))\} = \varphi\{\underset{n\to\infty}{Lt}\ (n(1 - g_n(t)))\} = \varphi\{\psi(t)\}.$$

This implies that $\omega(t) = e^{-\psi(t)}$ must be stable. Hence $f(t) = \varphi\{\psi(t)\}$ is $\varphi$-stable.

**Remark.4.1:** This extends Theorem.A.2 in Ramachandran (1997). One may also invoke Theorem.2.6 to prove this. Because by Ramachandran (1997) the condition

$$\underset{n\to\infty}{Lt}\ \{n(1- g(t/a_n)) + it\mu_n\} = \psi(t)\ \text{is equivalent to}\ \underset{n\to\infty}{Lt}\ \{n(1- g_n(t))\} = \psi(t).$$

**Theorem.4.2:** A CF $f(t)$ has a nonempty DP$\varphi$-A if and only if it is $\varphi$-ID.

**Proof:** Suppose $f(t)$ is $\varphi$-ID. Then there exists a CF $e^{-\psi(t)}$ that is ID such that $f(t) = \varphi\{\psi(t)\}$. Since $e^{-\psi(t)}$ is ID it has a DPA and let $g(t)$ belongs to it. Now by Lemma.2.2 and the transfer theorem,

$$f(t) = \varphi\{\psi(t)\} = \underset{k\to\infty}{Lt}\ \varphi\{n_k(1- g_k(t))\}.$$

Hence a CF that is $\varphi$-ID has a nonempty DP$\varphi$-A.

Conversely, let $f(t)$ has a nonempty DP$\varphi$-A. Then, as in (4.4) we have as $k \to \infty$, $\varphi\{n_k(1- g_k(t)) \to f(t)$.

Hence by the very definition $f(t)$ is $\varphi$-ID and the proof is complete. This result extends Theorems 4.1 and 4.4 in Mohan, *et al.* (1993).

**Theorem.4.3:** A CF $f(t)$ has a nonempty strict-sense DP$\varphi$-A where $\{n_k\}$ satisfies $\underset{k\to\infty}{Lt}\ \inf\ (n_k/n_{k+1}) >0$, if and only if it is $\varphi$-semi-stable.

**Proof:** This is another way of conceiving Theorem.2.7 and extends Theorem.2.5.1 in Sandhya (1991) (Theorem.3.1 in Sandhya and Pillai (1999)).

Another Look at Random Infinite Divisibility

These results also prove that a CF $g(t) \in$ DA (DPA) of the CF $e^{-\psi(t)}$ that is stable (semi-stable, ID) if and only if $g(t) \in$ D$\varphi$-A (DP$\varphi$-A) of $\varphi\{\psi(t)\}$ that is $\varphi$-stable ($\varphi$-semi-stable, $\varphi$-ID). Hence;

**Theorem.4.4:** The D$\varphi$-A of a $\varphi$-stable law with CF $f(t) = \varphi\{\psi(t)\}$ coincides with the DA of the stable law with CF $e^{-\psi(t)}$.

**Theorem.4.5:** The DP$\varphi$-A of a $\varphi$-ID law with CF $f(t) = \varphi\{\psi(t)\}$ coincides with the DPA of the ID law with CF $e^{-\psi(t)}$.

**Theorem.4.6:** The DP$\varphi$-A of a $\varphi$-semi-stable law with CF $f(t) = \varphi\{\psi(t)\}$ coincides with the DPA of the semi-stable law with CF $e^{-\psi(t)}$.

Also, since every stable law belongs to its own DA we can extend the analogous property of geometric stable laws *viz*. Theorem.2.4.1 in Sandhya (1991) (Theorem.2.2 in Sandhya and Pillai (1999)) and Theorem.A.1 in Ramachandran (1997) as:

**Theorem.4.7:** A $\varphi$-stable law belongs to its own D$\varphi$-A.

Next we have the apparently weaker requirement for $\varphi$-attraction analogous to Theorem.2.4.3 in Sandhya (1991) (Theorem.2.3 in Sandhya and Pillai (1999)) on geometric strictly stable laws.

**Theorem.4.8:** For a CF $g(t)$ to be $\varphi$-attracted to a $\varphi$-strictly stable law $f(t)$ it is necessary that $\varphi\{a_n(1 - g_k(t))\} \to f(t)$ and the integer part of $a_n$, $[a_n] = n$, $\forall n$. Conversely, if $\varphi\{a_n(1 - g_k(t))\} \to f(t)$ and $[a_n] = n$, $\forall n$, then $f(t)$ is $\varphi$-strictly stable and $g(t)$ is $\varphi$-attracted to $f(t)$.

**Remark 4.2:** When $\varphi$ is exponential we have the notions of exponential attraction and partial exponential attraction. But these are apparently different from the geometric attraction and partial geometric attraction developed in Sandhya (1991), Sandhya and Pillai (1999), Mohan, *et al.* (1993), Ramachandran (1997), GK (1996), and Klebanov and Rachev (1996). In these works the geometric attraction takes place when the geometric parameter $p$ runs through $1/n$ as $n$ increases to $\infty$ and partial geometric attraction when $p$ runs through a null sequence $\{p_n\}$. However the domain of geometric attraction (partial geometric attraction) coincides with D$\varphi$-A (DP$\varphi$-A) when $\varphi$ is exponential since both of them coincide with the DA (DPA) of the corresponding stable (ID) law. An extension of this argument shows that the DA (DPA) of a CF $f(t) = e^{-\psi(t)}$ coincides with the D$\varphi$-A (DP$\varphi$-A) of the CF $\varphi\{\psi(t)\}$ for any LT $\varphi$.

From the transfer theorem and Lemma.2.5 it is easy to see that for a given LT $\varphi$, $\varphi$-attraction and partial $\varphi$-attraction holds good with respect to every PGF in $\wp_\varphi$. Hence we may generalize these notions by specifying the requirement as:

**Definition.4.3:** For PGFs $P_\lambda \in \wp_\varphi$ with the distribution of the PGF $P_\lambda$ being independent of that of the CF $g(t)$; $g(t)$ belongs to the D$\varphi$-A of the CF $f(t)$ if;

$$\underset{n\to\infty}{Lt}\ P_n\{g_n(t)\} = f(t), \text{ for all } t \in \mathbf{R} \text{ (instead of (4.3)) and}$$

**Definition.4.4:** $g(t)$ belongs to the DP$\varphi$-A of $f(t)$ if;

$$\underset{k\to\infty}{Lt}\ P_{n_k}\{g_k(t)\} = f(t), \text{ for all } t \in \mathbf{R} \text{ (instead of (4.4)),}$$

where, $P_\lambda(s) = s^j \varphi\{\lambda(1-s^m)\}$, $0<s<1$, $\lambda \in (0,\infty)$ and $j \geq 0$ & $m \geq 1$ are integers.

However, whether they coincide with the DA or the DPA of the corresponding stable or ID laws (with CF $f(t) = e^{-\psi(t)}$) is not clear. The next result is in this direction where we establish a necessary and sufficient condition for the convergence to ID laws from that of an N-sum, that is an analogue of Theorem.4.6.5 in GK (1996, p.149) for PGFs in $\wp_\varphi$. From another angle it extends the formulation in Theorem.4.8 to PGFs in $\wp_\varphi$.

Let $\{X_{\theta,i}\}$ are i.i.d for every $\theta \in \Theta$ and $N_\theta$ be independent of $\{X_{\theta,i}\}$ for every $\theta \in \Theta$ with PGF $P_\theta \in \wp_\varphi$. Let $[1/\theta]$ denote the integer part of $1/\theta$. Then we have:

**Theorem.4.9:** Let $F(x)$ be a $\varphi$-ID d.f with CF $f(t) = \varphi(-m \log g(t))$, $m \geq 1$ integer. Then,

$$P\{\sum_{j=1}^{N_\theta} X_{\theta,j} < x\} \xrightarrow{w} F(x) \text{ as } \theta \to 0 \qquad (4.5)$$

if and only if there exists an ID d.f $G(x)$ with CF $g(t)$ and

$$P\{\sum_{j=1}^{[1/\theta]} X_{\theta,j} < x\} \xrightarrow{w} G(x) \text{ as } \theta \to 0, \qquad (4.6)$$

**Proof:** Our arguments follow that of GK (1996, p.149) except that leading to (4.9).

The sufficiency of the condition (4.6) follows from the transfer theorem by invoking the relation $\theta [1/\theta] \to 1$ as $\theta \to 0$ and the condition $\theta N_\theta \xrightarrow{L} mU$ proved in Lemma.2.5.

To prove the necessity of (4.6): Let $f_\theta(t)$ be the CF of $X_{\theta,i}$. Then in terms of CFs and PGFs (4.5) means

$$\underset{\theta \to 0}{Lt}\ P_\theta\{f_\theta(t)\} = f(t) = \varphi(-m \log g(t)). \qquad (4.7)$$

On the other hand, the CF of $\sum_{j=1}^{[1/\theta]} X_{\theta,j}$ is $\{f_\theta(t)\}^{[1/\theta]}$ and

$$\log \{f_\theta(t)\}^{[1/\theta]} = [1/\theta] \log\{1 - (1 - f_\theta(t))\}$$

$$= [1/\theta] (f_\theta(t) - 1) + \kappa [1/\theta] |f_\theta(t) - 1|^2, \text{ where } |\kappa| \leq 1. \qquad (4.8)$$

But (4.7) is possible only if $\underset{\theta \to 0}{Lt} [1 - \{f_\theta(t)\}^m] = 0$, $m \geq 1$ integer so that

$$\underset{\theta \to 0}{Lt}\, f_\theta(t) = 1 \quad \forall\, t \in \mathbf{R}. \tag{4.9}$$

Now, $\varphi\{(1-f_\theta(t))/\theta\}$ is a CF that is $\varphi$-ID for every $\theta \in \Theta$, and $\underset{\theta \to 0}{Lt}\, \varphi\{(1-f_\theta(t))/\theta\}$ is also $\varphi$-ID by Theorem.2.2. Hence there exists a CF $h(t)$ that is ID such that

$$\underset{\theta \to 0}{Lt}\, \{(1-f_\theta(t))/\theta\} = -\log h(t) \quad \forall\, t \in \mathbf{R}. \tag{4.10}$$

By (4.9) and (4.10) it follows from (4.8) that

$$\underset{\theta \to 0}{Lt}\, \{f_\theta(t)^{[1/\theta]}\} = h(t) \quad \forall\, t \in \mathbf{R}. \tag{4.11}$$

Applying the transfer theorem under conditions (4.11) and $\theta N_\theta \xrightarrow{L} mU$ proved in Lemma.2.5 it follows that the CF of $\sum_{j=1}^{N_\theta} X_{\theta,j}$ as $\theta \to 0$, is $\varphi\{-m \log h(t)\}$. Hence by (4.5) (or (4.7)) $h(t) \equiv g(t)$. That is, by (4.11), (4.6) is true with the d.f $G$ being ID, completing the proof.

**Remark.4.3** This theorem enables us to conclude that under the Definitions.4.3 & 4.4; if the CF $g(t)$ belongs to the D$\varphi$-A of a $\varphi$-stable law (DP$\varphi$-A of a $\varphi$-ID law) with CF $f(t) = \varphi\{\psi(t)\}$ then it is also a member of the DA of the stable law with CF $\omega(t) = e^{-\psi(t)}$ (DPA of the ID law with CF $\omega(t) = e^{-\psi(t)}$) and the converses are also true. All that we need is to prescribe, $[1/\theta] = n$ for $\varphi$-attraction and $[1/\theta] = n_k$ for partial $\varphi$-attraction. Hence if $g(t)$ is attracted to $\omega(t)$ then we have (4.6) holding good as $\theta \to 0$ through values such that $[1/\theta] = n$ and consequently (4.5) is true as $\theta$ runs through these values and hence $g(t)$ is $\varphi$-attracted to $f(t)$. Conversely, if (4.5) holds as $\theta$ runs through values such that $[1/\theta] = n$ then $g(t)$ is $\varphi$-attracted to $f(t)$, consequently (4.6) holds implying that $g(t)$ is attracted to $\omega(t)$. A similar line of argument holds true when $[1/\theta] = n_k$, a subsequence of the set of positive

Another Look at Random Infinite Divisibility

integers and we have the partial φ-attraction of $g(t)$ to $f(t)$ implying partial attraction of $g(t)$ to $\omega(t)$ and vice-versa. Thus the DA (DPA) coincides with the Dφ-A (DPφ-A) for each PGF $P_\theta \in \wp_\varphi$, and none of them are empty as well. Further, if the sequence $\{n_k\} = \{[1/\theta]\}$ satisfies $\underset{k \to \infty}{Lt} \inf (n_k/n_{k+1}) > 0$, then we have the DPφ-A of a φ-semi-stable law coinciding with the DPA of the corresponding semi-stable law. We can now formulate results that are 'φ-analogues' of the remarks in Section.4 in Mohan, *et al.* (1993) and Theorems.4.6.6 & 4.6.7 in GK (1996, p.151) for $P_\theta \in \wp_\varphi$. Satheesh and Sandhya (2003) has a result in this direction in terms of the CF $f_\theta(t)$ generalizing Theorem.1 in Feller (1971, p.555).

**Remark.4.4:** Thus we have extended and generalized the notions of ID laws, semi-stable laws, stable laws, attraction and partial attraction to N-sums for a class of PGFs for each given LT $\varphi$. This class of PGFs corresponding to $N$ is a much larger class than those under the assumption that the PGFs form a commutative semi-group.

**Acknowledgement:** Thanks are due to Dr. N Balakrishna of CUSAT for clarifying certain points in these developments in the seminars the author gave in the department in 2001. The author is thankful also to the referee for a careful reading, removing some mistakes and the valuable suggestions for improvement.

## References

1. **Bunge, J.** (**1996).** Composition semi groups and random stability, *Ann. Prob.*, 24, 1476 – 1489.

2. **Erdogan, M.B. and Ostrovskii, I.V.** (**1998).** Analytic and asymptotic properties of generalized Linnik probability densities, *J. Math. Anal. Appl.*, 217, 555-578.

3. **Feller, W. (1968).** *An Introduction to Probability Theory and Its Applications*, Vol.1, 3$^{rd}$ Edition, John Wiley and Sons, New York.

4. **Feller, W. (1971).** *An Introduction to Probability Theory and Its Applications*, Vol.2, 2$^{nd}$ Edition, John Wiley and Sons, New York.

5. **Gaver, D.P and Lewis, P.A.W. (1980).** First-order auto regressive gamma sequences and point processes, *Adv. Appl. Prob.*, 12, 727 – 745.

6. **Gnedenko, B.V. and Korelev, V.Yu. (1996).** *Random Summation, Limit Theorems and Applications*, CRC Press, Boca Raton.

7. **Grandell, J. (1997).** *Mixed Poisson Processes*, Chapman and Hall, London.

8. **Klebanov, L.B. and Rachev, S.T. (1996).** Sums of a random number of random variables and their approximations with ν-accompanying infinitely divisible laws, *Serdica Math. J.*, 22, 471 – 496.

9. **Klebanov, L.B; Maniya, G.M. and Melamed, I.A. (1984).** A problem of Zolotarev and analogues of infinitely divisible and stable distributions in the scheme of summing a random number of random variables, *Theor. Probab. Appl.*, 29, 791 – 794.

10. **Kozubowskii, T.J. and Panorska, A.K. (1996).** On moments and tail behavior of ν-stable random variables, *Statist. Prob. Letters*, 29, 307-315.

11. **Kozubowskii, T.J. and Panorska, A.K. (1998).** Weak limits for multivariate random sums, *J. Multi. Anal.*, 67, 398-413.

12. **Mohan, N.R; Vasudeva, R. and Hebbar, H.V. (1993).** On geometrically infinitely divisible laws and geometric domains of attraction, *Sankhya-A*, 55, 171 – 179.

13. **Ramachandran, B. (1997).** On geometric stable laws, a related property of stable processes, and stable densities of exponent one, *Ann. Inst. Statist. Math.*, 49, 299 – 313.

14. **Sandhya, E. (1991).** *Geometric Infinite Divisibility and Applications*, Ph.D. Thesis (unpublished), University of Kerala, January-1991.

15. **Sandhya, E. (1996).** On a generalization of geometric infinite divisibility, *Proc. $8^{th}$ Kerala Science Congress*, January-1996, 355 – 357.

16. **Sandhya, E. and Pillai, R.N. (1999).** On geometric infinite divisibility, *J. Kerala Statist. Assoc.*, 10, 1-7.

17. **Sandhya, E. and Satheesh, S. (1996).** On the membership of semi-$\alpha$-Laplace laws in class-L, *J. Ind. Statist. Assoc.*, 34, 77 – 78.

18. **Satheesh, S. (2002).** Aspects of Randomization in Infinitely Divisible and Max-Infinitely Divisible Laws, *ProbStat Models*, 1, 7 – 16.

19. **Satheesh, S. (2003).** A supplement to the Bose-Dasgupta-Rubin (2002) review of infinitely divisible laws and processes, *preprint*.

20. **Satheesh, S. and Nair, N.U. (2002).** Some classes of distributions on the non-negative lattice, *J. Ind. Statist. Assoc.*, 40, 41 – 58.

21. **Satheesh, S; Nair, N.U. and Sandhya, E. (2002).** Stability of random sums, *Stochastic Modeling and Applications*, 5, 17 –26.

22. **Satheesh, S. and Sandhya, E. (2003).** Infinite divisibility and max-infinite divisibility with random sample size, *Statistical Methods*, 5 (2), 126-139.

23. **Satheesh, S. and Sandhya, E. (2004).** Semi-stability of sums and maximums in samples of random size, *preprint*.

24. **Steutel, F.W. and van Harn, K. (1979).** Discrete analogues of self-decomposability and stability, *Ann. Prob.*, 7, 893-899.